\documentclass[12pt]{article}
\usepackage{amssymb}
\usepackage{amsmath}
\usepackage{amsfonts}
\usepackage{bbm}
\usepackage[mathscr]{euscript}
\ifx\pdftexversion\undefined
  \usepackage[dvips]{graphics}
\else
  \usepackage[pdftex]{graphics}
\fi

\usepackage{fullpage}

\newcommand{\R}{\mathbb{R}}
\newcommand{\E}{\mathbb{E}}
\newcommand{\N}{\mathbb{N}}
\newcommand{\Z}{\mathbb{Z}}
\newcommand{\p}{\mathbb{P}}

\newcommand{\ip}[2]{\left\langle{#1},{#2}\right\rangle}
\newcommand{\ket}[1]{|{#1}\rangle}
\newcommand{\bra}[1]{\langle{#1}|}

\usepackage{theorem}
\theorembodyfont{\upshape}
\newtheorem{thm}{Theorem}
\newtheorem{cor}{Corollary}

\renewcommand{\epsilon}{\varepsilon}

\title{A Thinning Analogue of de Finetti's Theorem}
\author{Shannon Starr
\footnote{\tt sstarr@math.mcgill.ca}\\
Centre de Recherches Math\'ematiques,
Universit\'e de Montr\'eal, and\\
Deparment of Mathematics and Statistics,
McGill University\\
Montr\'eal, Qu\'ebec, Canada}

\date{18 June 2004}

\begin{document}

\markright{Thinning analogue of de Finetti's theorem}

\maketitle

\begin{abstract}
We consider a notion of thinning for triangular arrays of random variables $(X^{(n)}_k\, :\, n \in \N_+,\, 1\leq k\leq n)$, taking values in
a compact metric space $\mathscr{X}$.
This is the sequence of transitions $T^{n+1}_n : \mathcal{M}^+_1(\mathscr{X}^{n+1}) \to \mathcal{M}^+_1(\mathscr{X}^{n})$
for each $n \in \N_+$, where $T^{n+1}_n \mu_{n+1}$ is the law of the random subsequence
$(X^{(n+1)}_{i_1},\dots,X^{(n+1)}_{i_n})$, where $1\leq i_1<\dots<i_n\leq n+1$ is chosen uniformly, at random.
We classify the set of all sequences of measures
$(\mu_n \in \mathcal{M}^+_1(\mathscr{X}^n)\, :\, n \in \N_+)$, which are thinning-invariant.
We prove that these form
a Bauer simplex and identify the extreme points.

While de Finetti's theorem is useful for identifying limit Gibbs states for classical
spin systems on complete graphs, our classification is useful for describing spin systems
on a one-dimensional chain, with a finite-ranged interaction, but where
the chain is normalized to have finite length.
I.e., $\Lambda_n = \{k/n\, :\, k=1,\dots,n\}$.
A similar statement is true when ``Gibbs state for classical spin system''
is replaced by ``invariant measure for asymmetric exclusion process''.

\end{abstract}

\tableofcontents

\section{Introduction}

Let us begin by considering de Finetti's theorem.

The usual statement of de Finetti's theorem is that a sequence of exchangeable
random variables is necessarily a mixture of independent, identically distributed
random variables.
Exchangeability means that, if the random variables $(X_n \, :\, n \in \N_+)$ take values
in $\mathscr{X}$ and $f : \mathscr{X}^k \to \R$ is any measurable function, for $k \in \N_+$,
then $\E\{f(X_{i_1},\dots,X_{i_k})\}$ does not depend on the values of the indices
$i_1,\dots,i_k \in \N_+$ as long as $i_1,\dots,i_k$ are distinct.
An excellent reference is \cite{HewittSavage}, wherein the theorem is proved by
showing that the extreme points must be i.i.d.\ sequences and using Choquet theory.
Another beautiful reference is \cite{Aldous}, which gives a deep picture of exchangeability,
as well as two different proofs (one of which is original to Aldous).

Now we will consider a different perspective.
This is motivated by Kingman's random partition structures,
\cite{Kingman}, but it is not necessary to be acquainted with those
to appreciate this perspective.

Consider a triangular sequence of random variables given
by a sequence of measures $(\mu_n \in \mathcal{M}^+_1(\mathscr{X}^n)\, :\, n\in \N_+)$.
I.e., $\mu_n$ is a probability measure describing the random variables
$(X^{(n)}_1,\dots,X^{(n)}_n) \in \mathscr{X}^n$.
For $k,n \in \N_+$ with $k\leq n$, consider the transition
$\widetilde{T}^n_k : \mathcal{M}^+_1(\mathscr{X}^n) \to \mathcal{M}^+_1(\mathscr{X}^k)$,
which is uniquely defined by specifying
\begin{multline}
  (\widetilde{T}^{n}_k \mu_n)(\{(X_1 \in A_1,\ldots, X_k \in A_k\})\, =\\ 
  \binom{n}{k}^{-1} \sum_{1\leq i_1<\dots<i_k\leq n} \frac{1}{k!} \sum_{\pi \in S_k}
  \mu(\{X_{i_{\pi(1)}} \in A_1,\ldots,X_{i_{\pi(k)}} \in A_k\})
\end{multline}
for all measurable events $A_1,\dots,A_k$, subsets of $\mathscr{X}$.
In other words, if the random variables  $(X_1,\dots,X_n)$,
are distributed according to $\mu_n$, then the law of the random sequence
$(\widehat{X}_1,\dots,\widehat{X}_k)$ is 
given by $(\widetilde{T}^n_k \mu_n)$,
if $(\widehat{X}_1,\dots,\widehat{X}_k)$  is obtained
by selecting a uniform, random $k$-subset $(X_{i_1},\dots,X_{i_k})$,
and then shuffling by a uniform, random permutation $\pi^{-1} \in S_k$.

We define the sequence $(\mu_n \in \mathcal{M}_1(\mathscr{X}^n)\, :\, n \in \N_+)$
to describe an ``exchangeable triangular array'' of random variables if
$\widetilde{T}^n_k \mu_n = \mu_k$ for all $k,n \in \N_+$ with $k\leq n$.
In particular, note that in this case each $\mu_n$ is exchangeable in the sense that it is 
invariant with respect to the natural action of $S_n$.

In this case, similarly to the case above, one can prove that each exchangeable 
triangular array $(\mu_n\, :\, n \in \N_+)$ can be uniquely represented
as a mixture of triangular arrays of i.i.d.\ random variables:
\begin{multline}
\label{def:eq}
  \exists\, ! P \in \mathcal{M}^+_1(\mathcal{M}^+_1(\mathscr{X})) \mbox{ s.t. }
  \forall n \in \N_+\, ,\\
  \mu_n(dX_1\otimes \cdots\otimes dX_n)\, =\, \int_{\mathcal{M}^+_1(\mathscr{X})} 
  \left[\prod_{k=1}^n \mu(dX_k)\right]\, P(d\mu)\, .
\end{multline}
This holds in the same generality as the result in \cite{HewittSavage},
and could be proved by the same technique.

This version of de Finetti's theorem suggests the question:
``What if, in the definition of $\widetilde{T}^n_k$, we
did not shuffle by a uniform, random permutation $\pi \in S_k$?''
A more general question is to suppose that $\pi$ is chosen by another
distribution than either $\delta_e$ or the Haar measure on $S_n$.
We will not consider the more general problem, because we believe that under
reasonable conditions one will obtain exchangeable distributions,
unless one takes $\delta_e$ as the measure on $S_n$.

Let us define the transition
$T^n_k : \mathcal{M}^+_1(\mathscr{X}^n) \to \mathcal{M}^+_1(\mathscr{X}^k)$,
such that
\begin{multline}
  (T^{n}_k \mu_n)(\{X_1 \in A_1,\ldots, X_k \in A_k\})\, =\\ 
  \binom{n}{k}^{-1} \sum_{1\leq i_1<\dots<i_k\leq n} 
  \mu(\{X_{i_1} \in A_1,\ldots,X_{i_k} \in A_k\})
\end{multline}
for all measurable events $A_1,\dots,A_k \subset \mathscr{X}$.
We call a sequence of measures, $(\mu_n \in \mathcal{M}^+_1(\mathscr{X}^n)\, :\, n\in \N_+)$,
a thinning invariant triangular array if, for all $k,n \in \N_+$ with $k\leq n$,
it is true $\mu_k = T^{n}_k \mu_n$.

Of course, every exchangeable triangular array is also thinning invariant,
but the reverse is not true.

\section{Main Result}

The main result is quite simple to understand, and is a very natural extension 
of the conclusion of de Finetti's theorem.
However, whereas in understanding de Finetti's theorem, the basic structures
to consider are measures $\mu \in \mathcal{M}^+_1(\mathscr{X})$
(as well as the mixtures of those, $\mathcal{M}^+_1(\mathcal{M}^+_1(\mathscr{X}))$)
in the present context, the basic structures are Borel measurable functions
from $[0,1]$ into $\mathcal{M}^+_1(\mathscr{X})$, as we will see.
Therefore, we will try to be careful when it comes to topology and measure theory.

We consider $[0,1]$ as a metric space, and equip it with standard Lebesgue measure.
We will use the Borel $\sigma$-algebras for $[0,1]$ as well as for $\mathscr{X}$.

Note that since $\mathscr{X}$ is compact, it is obviously separable.
By the Arzela-Ascoli theorem, $\mathcal{C}(\mathscr{X})$ is locally
compact, hence it is also separable.
This will be useful later.
Obviously $\mathcal{M}^+_1(\mathscr{X})$
is convex. 
It is compact, in the vague topology, by the Banach-Alaoglu theorem.
In fact, this topology is also metrizable, using the Riesz-Markov theorem
and the fact that $\mathcal{C}(\mathscr{X})$ is separable.
This fact will also be useful later.

The set $\mathcal{B}([0,1],\mathcal{C}(\mathscr{X}))$ is defined to be
the bounded, Borel measurable functions from $[0,1]$ to $\mathcal{C}(\mathscr{X})$,
defined for Lebesgue-a.e.\ $t \in [0,1]$.
If we write such a function as $F(t,x)$ for $t \in [0,1]$, $x \in \mathscr{X}$,
then to be Borel measurable means that for any $f \in \mathcal{C}(\mathscr{X})$
and any $\epsilon>0$, the subset of $[0,1]$ defined as
$\{t\, :\, \sup_{x \in \mathscr{X}} |F(t,x)-f(x)| < \epsilon\}$ is a Borel set.
We define
\begin{equation}
  \mathcal{L}\, =\, L^1([0,1],\mathcal{C}(\mathscr{X}))\, ,
\end{equation}
the Banach space whose norm is
\begin{equation}
  \|F\|_{\mathcal{L}}\, =\, \int_0^1 \max_{x \in \mathscr{X}} |F(t,x)|\, dt\, .
\end{equation}
More precisely, $\mathcal{L}$ is defined as the norm-closure,
in the set of functions
$\mathcal{B}([0,1],\mathcal{C}(\mathscr{X}))$, of all functions
$\sum_{n=1}^N \mathbbm{1}_{E_n}(t) f_n(x)$, where $N$ is any 
positive integer $N \in \N_+$, and  
$E_1,\dots,E_N$ are Borel measurable subsets of $[0,1]$, 
and $f_1,\dots,f_N$ are continuous functions in $\mathcal{C}(\mathscr{X})$.

As noted above, $\mathcal{C}(\mathscr{X})$ is separable.
Also note that there is a is a countable collection of Borel measurable subsets of 
$[0,1]$ (e.g., open intervals with rational endpoints),
such that any Borel measurable set $E$ can be arbitrarily well-approximated
in Lebesgue-measure by a sequence of subsets $E_n$, such that each $E_n$ is a finite
union of these countable number of sets.
Therefore, from its definition as the closure of simple functions, 
it is clear that $\mathcal{L}$ is separable.

We define a convex set $\mathcal{K}$ to be the set
of Borel probability measures $\nu \in \mathcal{M}^+_1(\mathscr{X}\times [0,1])$,
satisfying that for every Borel measurable set $E \subset [0,1]$,
it is true that 
\begin{equation}
  \nu(\mathscr{X}\times E)\, =\, |E|\, =\, \int_0^1 \mathbbm{1}_E(t)\, dt\, .
\end{equation}
(We denote Lebesgue measure of $E$ by $|E|$, not $m(E)$.)

Given any $f \in \mathcal{C}(\mathscr{X})$, it is clear that
\begin{equation}
  E \mapsto \int_{\mathscr{X} \times [0,1]} f(x) \mathbbm{1}_E(t)\, \nu(dx\otimes dt)
\end{equation}
defines a signed measure, which is absolutely continuous with respect to
Lebesgue measure.
We denote the density with respect to Lebesgue measure as 
$\ell(f,\cdot) \in L^1(dt)$, so that
\begin{equation}
  \int_{\mathscr{X} \times [0,1]} f(x) \mathbbm{1}_E(t)\, \nu(dx\otimes dt)\,
  =\, \int_0^1 \ell(f,t)\, \mathbbm{1}_E(t)\, dt\, .
\end{equation}
If, additionally, $f>0$, then $\ell(f,t)>0$ for Lebesgue-a.e.\ $t \in [0,1]$.
For every Borel measurable $E \subset [0,1]$, one has
\begin{equation}
  \left|\int_0^1 \ell(f,t)\, \mathbbm{1}_E(t)\, dt\right| \leq \|f\|_{\mathcal{C}(\mathscr{X})} \cdot |E|\, ,
\end{equation}
using H\"older's inequality.
Therefore, by density, we conclude that $|\ell(f,t)| \leq \|f\|$ for a.e.\ $t \in [0,1]$.
A similar argument shows that $\ell(f,t)$ is linear in $f$ for almost every
$t$.
Note that, since $\mathcal{C}(\mathscr{X})$ is separable and $\ell(\cdot,t)$ is
uniformly continuous, linearity in $f$ can be made to hold simultaneously for 
all of $\mathcal{C}(\mathscr{X})$, for a.e.\ $t \in [0,1]$,
because a countable union of null sets is null.
Therefore, by the Riesz-Markov theorem, for almost every $t \in [0,1]$,
we can identify a measure $\mu(\cdot,t) \in \mathcal{M}^+_1(\mathscr{X})$ such that
$\nu(dx\otimes dt) = \mu(dx,t)\, dt$.
Note that, since the density $\ell(f,t)$ is chosen as a Borel measurable function
of $t$,
it is clear that for any $f \in \mathcal{C}(\mathscr{X})$, the function
\begin{equation}
\label{bor:eq}
  t \mapsto \int_{\mathscr{X}} f(x)\, \mu(dx,t)
\end{equation}
is Borel measurable.

Henceforth we will think of $\mathcal{K}$ as the set of
Borel measurable functions from $[0,1] \to \mathcal{M}^+_1(\mathscr{X})$,
defined for Lebesgue-a.e.\ $t \in [0,1]$.
We write this as $\mu(dx,t)$, as above.
Note that since the topology of $\mathcal{M}^+_1(\mathscr{X})$ is the vague
topology, to be Borel means precisely that the map 
defined in (\ref{bor:eq}) is Borel, for every $f \in \mathcal{C}(\mathscr{X})$.
Next we will endow $\mathcal{K}$ with a topology.

Given any $F \in \mathcal{L}$ and any $\mu \in \mathcal{K}$, we define the number
\begin{equation}
  \ip{\mu}{F}\, :=\, \int_0^1 \int_{\mathscr{X}} F(x,t)\, \mu(dx,t)\, dt\, .
\end{equation}
This is a bilinear map on $\mathcal{K} \otimes \mathcal{L}$.
It is continuous as a function on $\mathcal{L}$: 
in fact, since $\mu(\cdot,t) \in \mathcal{M}^+_1(\mathscr{X})$
for almost every $t \in [0,1]$, it is clear that $|\ip{\mu}{F}| \leq \|F\|_{\mathcal{L}}$.
I.e., $\ip{\mu}{\cdot}$ is a bounded, linear functional on $\mathcal{L}$.
We define the topology of $\mathcal{K}$ to be the weak topology
with respect to all the maps $\ip{\cdot}{F}$, for $F \in \mathcal{L}$.
It will be important to us that, with this topology, $\mathcal{K}$
is compact and metrizable.

The reason that $\mathcal{K}$ is metrizable is just that $\mathcal{L}$ is separable, as we noted
above.
Take $(F_n\, :\, n \in \N_+)$ to be a countable dense subset of the unit
ball in $\mathcal{L}$.
Then we can define a metric
\begin{equation}
  d(\mu,\nu)\, =\, \sum_{n=1}^\infty 2^{-n}\, |\ip{\mu}{F_n} - \ip{\nu}{F_n}|\, .
\end{equation}
A sequence $(\mu_n\, :\, n \in \N_+)$ converges to $\mu$ in $\mathcal{K}$
iff $\lim_{n \to \infty} \ip{\mu_n}{F} = \ip{\mu}{F}$ for every $F \in \mathcal{L}$,
and this happens iff $\lim_{n \to \infty} d(\mu_n,\mu) = 0$.
The fact that $\mathcal{K}$ is compact then follows by the same argument
as the Banach-Alaoglu theorem.

We are now ready to consider the fundamental example of a thinning-invariant
sequence of measures $(\mu_n \in \mathcal{M}^+_1(\mathscr{X}^n)\, :\, n \in \N_+)$.
Let $\mu \in \mathcal{K}$ be arbitrary.
For each $n \in \N_+$, we define a measure $\mu_n$ as
\begin{equation}
  \mu_n(dX_1\otimes\cdots\otimes dX_n)\, 
  =\,\int_{\Delta^{n}} \mu(dX_1,t_1)\cdots \mu(dX_n,t_n)\, 
  n!\, dt_1\cdots dt_n\, ,
\end{equation}
where $\Delta^{n}=\{(t_1,\dots,t_n)\, :\, 0<t_1<\dots<t_n<1\}$,
which has Lebesgue measure $1/n!$.
This is clearly a well-defined measure, because for any $f\in \mathcal{C}(\mathscr{X}^n)$
we know that the the following function
is Borel measurable and bounded on $[0,1]^n$:
\begin{equation}
  (t_1,\dots,t_n) \mapsto \int_{\mathscr{X}^n} f(X_1,\dots,X_n)\, \mu(dX_1,t_1)\cdots \mu(dX_n,t_n)\, .
\end{equation}
Therefore the following integral is well-defined:
\begin{multline}
  \int_{\mathscr{X}^n} f(X_1,\dots,X_n)\,  \mu_n(dX_1\otimes\cdots\otimes dX_n)\, 
  =\\ \int_{[0,1]^n} 
  \left[\int_{\mathscr{X}^n} f(X_1,\dots,X_n)\, \mu(dX_1,t_1)\cdots \mu(dX_n,t_n)\right]\\
  n!\, \mathbbm{1}[0<t_1<\cdots<t_n<1]\, dt_1\cdots dt_n\, .
\end{multline}
Hence, $\mu_n \in \mathcal{M}^+_1(\mathscr{X}^n)$ can be defined by the Riesz-Markov theorem.

The probabilistic interpretation is that $t_1,\dots,t_n$ are random variables distributed
by choosing $n$ points independently and uniformly in $[0,1]$ and then sorting
in order from least to greatest.
Conditional on the values of the $t_1<\dots<t_n$, one considers $X_1,\dots,X_n$
to be distributed according to the conditionally independent (but not identically
distributed) measure $\mu(dX_1,t_1)\cdots\mu(dX_n,t_n)$.
Then $\mu_n$ is the distribution of $(X_1,\dots,X_n)$,
averaged over the distribution of $(t_1,\dots,t_n)$.

If one defines $\rho_n(dt_1\otimes\cdots\otimes dt_n)$ to be the distribution of $(t_1,\dots,t_n)$,
as described above,
then clearly $(\rho_n \in \mathcal{M}^+_1([0,1]^n)\, :\, n \in \N_+)$ is thinning-invariant.
This is sufficient to prove that $(\mu_n\, :\, n \in \N_+)$ is thinning-invariant
by construction.

\begin{thm}
Let $(\mu_n\, :\, n \in \N_+)$ describe a thinning-invariant, triangular array of 
random variables.
Then there exists a unique Borel probability measure $P \in \mathcal{M}^+_1(\mathcal{K})$
such that for each $n \in \N_+$,
\begin{multline}
\label{thmeq}
  \mu_n(dX_1\otimes\cdots\otimes dX_n)\\
  =\,
  \int_{\mathcal{K}} \int_{\Delta^{n}} 
  \mu(dX_1,t_1)\cdots \mu(dX_n,t_n)\, 
  n!\, dt_1\cdots dt_n\, P(d\mu)\, .
\end{multline}
\end{thm}

{\it Remark 1.}
We can deduce the second version of de Finetti's theorem from this 
theorem.
Suppose that $(\mu_n\, :\, n \in \N_+)$ is exchangeable.
Then, in particular it is thinning-invariant.
So there is a unique Borel measure  $P \in \mathcal{M}_+(\mathcal{K})$
such that (\ref{thmeq}) holds.
But because of exchangeability, we can average over the $S_n$ orbit
of $\Delta^{n}$, instead,
which is just $[0,1]^n$.
Therefore, it is clear that $\mu_n$ is a mixture of product measures,
each factor just being the average of $\int_0^1 \mu(dx,t)\, dt$.
Because of uniqueness of $P$, it is clear, a posteriori, that
$P$ must have been concentrated on Lebesgue-a.e.\ constant functions
$\mu(dx,t) \equiv \mu(dx)$.

{\it Remark 2.} 
If one considers the simple example of $\rho_n$, then this is an extreme
point, corresponding to the measure-valued function $\mu(dx,t) = \delta(x-t)$.
Note that this is a continuous function of $t$,
with respect to the vague topology (although nowhere continuous
with respect to the total-variation norm).

\section{Proof}

The proof follows a well-known and standard
technique in probabilistic results of this kind.
For example, it is the technique used by Kingman
to prove his representation theorem, in \cite{Kingman}.
The idea could
be described rather briefly, as it is done, there.
But we will attempt to explain every detail.

\subsection{Probabilistic tools}

There are two basic tools which we use, which are the 
Kolmogorov extension theorem, and Doob's
reversed martingale convergence theorem.

The Kolmogorov extension theorem says that
a projective family of probability measures has a unique
projective limit, under certain conditions.
We will use a variant of the ordinary version.
See, for example, \cite{Stroock}, Exercise 3.1.18(iii) for a precise 
statement of the usual Kolmogorov extension principle,
which is the reference we follow.

To be projective essentially means that there is a family of 
probability measures
$(\mu_n \in \mathcal{M}^+_1(\mathcal{F}_n)\, :\, n \in \N_+)$,
and a semigroup of transitions
\begin{equation}
  \phi^n_k\, :\, \mathcal{M}^+_1(\mathcal{F}_n) \to \mathcal{M}^+_1(\mathcal{F}_k)\, ,
\end{equation}
such that $\phi^n_k \mu_n = \mu_k$ for all $1\leq k\leq n$.
We will consider the case that $\mathcal{F}_n$ is the Borel $\sigma$-algebra
for a product $\prod_{k=1}^n \mathscr{X}_k$, for some compact metric spaces
$(\mathscr{X}_k\, :\, k \in \N_+)$.
In this case the map $\phi^n_k$ is the projection,
\begin{equation}
  (\phi^n_k \mu_n)(A) = \mu_n\left(A \times \prod_{j=k+1}^n \mathscr{X}_j\right)\, .
\end{equation}
In \cite{Stroock}, the reader is led to a proof of the extension principle
when all $\mathscr{X}_n = \mathscr{X}$ for a single $\mathscr{X}$.

Let us give a concrete method for generalizing from that case 
to the case where the $\mathscr{X}_n$'s are allowed to vary.
Define $\widehat{\mathscr{X}} = \prod_{n=1}^\infty \mathscr{X}_n$,
with the product topology, which is compact by Tychonov's theorem.
Choose an $e_n \in \mathscr{X}_n$ for each $n \in \N_+$.
Given $\mu_n \in \bigotimes_{k=1}^n \mathcal{M}^+_1(\mathscr{X}_k)$, 
define $\widehat{\mu}_n \in \mathcal{M}^+_1(\widehat{\mathscr{X}}^n)$
such that, given measurable sets $A_k \subset \mathscr{X}_k$
and $B_k \in \prod_{j\neq k} \mathscr{X}_j$ for each $k=1,\dots,n$, we have
\begin{equation}
  \widehat{\mu}_n\left(\prod_{k=1}^n (A_k \times B_k)\right)\, 
  =\, \mu_n\left(\prod_{k=1}^n A_k\right) \times \prod_{k=1}^n \left(\bigotimes_{j\neq k} \delta_{e_j}\right)(B_k)\, .
\end{equation}
If $(\mu_n\, :\, n \in \N_+)$ is a projective system, then 
$(\widehat{\mu}_n\, :\, n \in \N_+)$ is as well, and of the form
considered in \cite{Stroock}.
Restricting attention to the $\sigma$-algebra generated by
cylinder sets of the form $\prod_{k=1}^n (A_k \times B_k)$ with $A_k \subset \mathscr{X}_k$
and $B_k = \prod_{j\neq k} \mathscr{X}_j$, then gives us the result
we desire.

It is quite likely that there is a purely topological argument, based on projective
systems, that bypasses all of these technicalities, but we did not find
it in the elementary textbooks we consulted.

The second tool from probability which we need is Doob's reversed martingale
theorem.
Again, the reader can consult Stroock's text, \cite{Stroock}, Exercise 5.2.40(iii).
Given a family of decreasing $\sigma$-algebra, $(\mathcal{F}_n\, :\, n \in \N_+)$,
and a Borel probability measure $\mu$ on $\mathcal{F}_1$, a 
sequence of random variables $(X_n \, :\, n \in \N_+)$ is a reversed martingale
(with respect to $\mu$)
if $X_n$ is $\mathcal{F}_n$-measurable and in fact $X_n \in L^1(\mathcal{F}_n,d\mu)$
for each $n$, and if
\begin{equation}
  \E^{\mu}\{X_n\, |\, \mathcal{F}_N\} = X_N
\end{equation}
for all $1\leq n\leq N$.
Note that this is the opposite of a normal martingale, hence the term
reversed martingale. (Of course Doob also proved a convergence theorem
for normal martingales.)

Doob's convergence theorem says that there is a $\mu$-almost surely
unique limit $X$, which is 
measurable with respect to $\mathcal{F}_\infty = \cap_{n=1}^\infty \mathcal{F}_n$,
such that $X_n$ converges to $X$, pointwise, $\mu$-a.s., and that $X \in L^1(\mathcal{F}_\infty,d\mu)$.
In our case, we will apply this when each $X_k$ is uniformly bounded, all by
the same bound, in which case it is clear that the limit also satisfies
the same bound, almost surely.

Finally, we observe that, since a countable union of $\mu$-null sets is 
$\mu$-null,
we can apply the convergence theorem simultaneously to a countable
union of reversed martingales.
We will need this generalization, when we construct a function from its moments.

\subsection{Application of probabilistic tools}

Suppose that $(\mu_n\, :\, n \in \N_+)$ is thinning-invariant.
Then, for any $N \in \N_+$, we could define a measure 
$\mu^{(N)} \in \bigotimes_{n=1}^N \mathcal{M}^+_1(\mathscr{X}^n)$
on the triangular sequence of random variables $(X^{(n)}_k\, :\, 1\leq n\leq N,\, 1\leq k\leq n)$,
as follows.

Let $(X^{(N)}_1,\dots,X^{(N)}_N)$ be distributed according to $\mu_N$.
For each $n$ such that $1\leq n\leq N-1$, conditional on the values of $(X^{(n+1)}_1,\dots,X^{(n+1)}_{n+1})$,
let $(X^{(n)}_1,\dots,X^{(n)}_n)$ be distributed as $(X^{(n+1)}_{i_1},\dots,X^{(n+1)}_{i_n})$
where $1\leq i_1<\dots<i_n\leq n+1$ is a uniformly chosen $n$-subset
of $\{1,\dots,n+1\}$.
Then, because of thinning-invariance, if we restrict attention to 
$(X^{k}_j\, :\, 1\leq k\leq n,\, 1\leq j\leq k)$ (with $n\leq N$)
the marginal of $\mu^{(N)}$ is equal to $\mu^{(n)}$.
Thus the sequence $(\mu^{(N)}\, :\, N \in \N_+)$ forms a projective
sequence of probability measures.
By the Kolmogorov extension theorem, 
there is a unique extension,
which we will call $\mu^{(\infty)}$.
It is a probability measure in
$\bigotimes_{n=1}^\infty \mathcal{M}^+_1(\mathscr{X}^n)$.

Define $\mathcal{F}_n$ to be the $\sigma$-algebra generated by
the random variables 
\begin{equation}
\mathcal{F}_n\, =\, \sigma(X^{(N)}_k\, :\, N\geq n,\, 1\leq k\leq N)\, 
\end{equation}
for each $n \in \N_+$.
So $(\mathcal{F}_n\, :\, n \in \N_+)$ is a decreasing family of
$\sigma$-algebras, and $\mu^{(\infty)}$ is a measure
on $\mathcal{F}_1$.
Define the tail $\sigma$-field as 
\begin{equation}
  \mathcal{F}_\infty\, =\, \bigcap_{N=1}^\infty \mathcal{F}_N\, .
\end{equation}
Let $n \in \N_+$ and let $f \in \mathcal{B}(\mathscr{X}^n,\R)$
be a bounded, Borel measurable function on $\mathscr{X}^n$.
Then we can define a reversed martingale which we call
$(\phi_N[f]\, :\, N\geq n)$ by
\begin{equation}
  \phi_N[f] = \E\{f(X_1,\dots,X_n)\, |\, \mathcal{F}_N\}\, .
\end{equation}

It is clear that $\phi_N[f]$ is bounded.
Therefore, by Doob's reversed martingale convergence theorem,
there exists a random variable, $\phi_\infty[f]$, measurable with
respect to the tail $\sigma$-field $\mathcal{F}_\infty$,
such that 
\begin{equation}
  \lim_{N \to \infty} \phi_N[f]\, =\, \phi_\infty[f]\, ,\quad \mu^{(\infty)}-\mbox{a.s.}
\end{equation}
We will use this fact to construct the measure-valued function, $\mu(dx,t)$,
which will be measurable with respect to $\mathcal{F}_\infty$.

\subsection{Construction of the measure on finite algebras}

Suppose that $R \in \N_+$ and that $(A_1,\dots,A_R)$ is a Borel measurable
partition of $\mathscr{X}$.
I.e., that $A_1,\dots,A_R$ are Borel measurable subset of $\mathscr{X}$, and 
\begin{equation}
  \forall r\neq s,\ A_r \cap A_s = \emptyset\quad \mbox{ and }\quad
  \bigcup_{r=1}^R A_r = \mathscr{X}\, .
\end{equation}
We will start by constructing the projection of $\mu(dx,t)$, restricted to the
$\sigma$-algebra generated by this partition.

Let us define a vector-valued function
\begin{equation}
  m(x)\, =\, (\mathbbm{1}[x \in A_1],\dots,\mathbbm{1}[x \in A_R])\, .
\end{equation}
Of course this function is Borel measurable and bounded, in fact its
range is the set of extreme points of the simplex 
\begin{equation}
  \nabla^{R-1}\, =\, \{v \in \R^R\, :\, v_r \geq 0 \mbox{ for } r=1,\dots,R,
  \mbox{ and } \sum_{r=1}^R v_r = 1\}\, .
\end{equation}
Note that this simplex is naturally isomorphic to the set of measures on
the $\sigma$-algebra spanned by the partition, $\sigma(A_1,\dots,A_R)$.
Namely, given $v \in \nabla^{R-1}$, we can define a measure
\begin{equation}
  \mu_v(E) = \sum_{r=1}^R v_r \mathbbm{1}[A_r \subset E]\, .
\end{equation}
For each $n\in \N_+$, let us define a random, vector-valued function
$M_n\, :\, [0,1) \to \nabla^{R-1}$, given by 
\begin{equation}
  M_n(t)\, =\, \sum_{k=1}^n m(X^{(n)}_k)\, \mathbbm{1}\hspace{-2pt} \left[\frac{k-1}{n}\leq t<\frac{k}{n}\right]\, .
\end{equation}
We will show that the functions $M_n$ converge in 
the weak $L^1([0,1],\R^R)$-topology, which is the closest
analogue of the weak-$\mathcal{L}$ topology we can get.

To do so, we consider some $\mathcal{F}_\infty$-measurable
statistics of the random triangular array
$\{X^{(n)}_k\, :\, n\in \N_+,\, 1\leq k\leq n\}$.
Namely, we define, for each $1\leq k\leq n$, 
\begin{equation}
  \beta(k,n-k+1)\, :=\, \phi_\infty[m(X^{(n)}_k)]\, .
\end{equation}
So $\beta(k,n-k+1)$ is a $\nabla^{R-1}$-valued random variable, by convexity, 
and it is measurable with respect to
$\mathcal{F}_\infty$.
Since $(X^{(n)}_1,\dots,X^{(n)}_n)$ is obtained from
$(X^{(N)}_1,\dots,X^{(N)}_N)$ 
by choosing a random $k$ subset,
we know that
\begin{equation}
  \phi_N[m(X^{(n)}_k)]\, =\, \sum_{K=1}^N p(N,K;n,k)\, m(X^{(N)}_K)\, ,
\end{equation}
where
\begin{equation}
  p(N,K;n,k)\, =\, \binom{K-1}{k-1}\, \binom{N-K}{n-k} \Bigg/ \binom{N}{n}\, .
\end{equation}

Let us define the Beta-kernel, $\Phi_{k,n-k+1} : [0,1] \to \R_+$, as 
\begin{equation}
  \Phi_{k,n-k+1}(t)\, =\, B(k,n-k+1)^{-1}\, t^{k-1}\, (1-t)^{n-k}\, 
\end{equation}
where $B(k,n-k+1)$ is the Beta-integral
\begin{equation}
  B(k,n-k+1)\, =\, \int_0^1 t^{k-1} (1-t)^{n-k}\, dt\, =\, (k-1)!\, (n-k)! / n!\, .
\end{equation}
Then it is easy to see that 
\begin{equation}
  \lim_{N \to \infty} \sum_{K=1}^\infty |p(N,K;n,k)\, -\, \Phi_{k,n-k+1}(K/N)|\, =\, 0\, .
\end{equation}
So, by H\"older's inequality,
\begin{equation}
  \lim_{N \to \infty} \left|\phi_N[m(X^{(n)}_k)]\, 
  - \int_0^1 M_N(t)\, \Phi_{k,n-k+1}(t)\, dt\right|\, =\, 0\, .
\end{equation}
In other words, since $\phi_N[m(X^{(n)}_k)]$ converges,
\begin{equation}
  \lim_{N \to \infty} \int_0^1 M_N(t)\, \Phi_{k,n-k+1}(t)\, dt\, =\, \beta(k,n-k+1)\, ,\quad 
  \mu^{(\infty)}-\mbox{a.s.}
\end{equation}
Setting $k=n$, we obtain, for each $n \in \N_+$,
\begin{equation}
  \lim_{N \to \infty}
  \int_0^1 M_N(t)\, t^{n-1}\, dt\, =\, \beta(n,1)/n\, ,\quad 
  \mu^{(\infty)}-\mbox{a.s.}
\end{equation}

We will use orthogonal polynomials to show that the sequence
$M_N$ converges in the weak $L^2([0,1],\R^R)$ topology.
From the above, it is obvious that for any polynomial $\varphi : [0,1] \to \R$,
the vectors $\ip{M_N}{\varphi}$ converge, and the limit is a 
linear combination of the vectors $\beta(n,1)$ for $1\leq n\leq \deg(\varphi)+1$.
Let $(\varphi_j\, :\, j \in \N)$ be the orthonormal polynomials 
for $[0,1]$ with Lebesgue measure.
By the Parseval relation, for each $N \in \N_+$.
\begin{equation}
  \sum_{j=0}^\infty \|\ip{M_N}{\varphi_j}\|^2\, =\, \|M_N\|_{L^2}\, =\, 1\, .
\end{equation}
Therefore, by Fatou's lemma, we have that
\begin{equation}
  \sum_{j=0}^\infty \|\lim_{N \to \infty} \ip{M_N}{\varphi_j}\|^2\, \leq\, 1\, .
\end{equation}
(See (\ref{defnorm}) for the precise norm.)
Hence, we can define an $L^2$-function $M(t)$ with these components
in the basis  $(\varphi_j\, :\, j \in \N)$.
By construction, $M_n$ converge to $M$, in the weak $L^2([0,1],\R^R)$
topology.

Next, we want to check that $M$ takes values in $\nabla^{R-1}$,
almost surely, which is not obvious since $M$ has been
constructed in terms of its moments.
But for any Borel measurable $E \subset [0,1]$, 
and any $r\in\{1,\dots,R\}$, the function 
$\varphi(t) = \widehat{e}_r \mathbbm{1}_E(t)$ is in $L^2([0,1],\R^R)$.
Therefore, since it is evidently true that
$\ip{M_N}{\varphi} \in [0,|E|]$, the same is true, $\ip{M}{\varphi} \in [0,|E|]$.
Running over all Borel measurable sets $E$, it is then apparent that
$(M(t),\widehat{e}_r) \in [0,1]$ for almost all $t$, $\mu^{\infty}$-a.s.
Similarly, one can show that since
\begin{equation}
  (M_N(t),\widehat{e}_1+\cdots+\widehat{e}_R)\, =\, 1\, ,\quad \mbox{for a.e. } t\in[0,1]\, ,
\end{equation}
that the same is true of $M$, $\mu^{\infty}$-a.s.

Finally, to show that $M_N \to M$ in weak $L^1([0,1],\R^R)$,
it suffices to observe that for any Borel measurable function
$F : [0,1] \to \nabla^{R-1}$, and any two $f,g \in L^1([0,1],\R^R)$,
it is true that $|\ip{f}{F} - \ip{g}{F}| \leq \|f-g\|_1$, where
we define 
\begin{equation}
\label{defnorm}
  \|f\|_p\, =\, \left[\int_0^1 \|f(t)\|_{p'}^{p}\, dt\right]^{1/p}\, ,
\end{equation}
for $p=1,2$, where $(1/p)+(1/p')=1$.
Then, by density of $L^2([0,1],\R^R)$ in $L^1([0,1],\R^R)$,
it follows that $M_N \to M$ in weak $L^1$
because it converges in weak $L^2$.

We can view $M$ as a function from $[0,1]$ to $\mathcal{M}^+_1(\sigma(A_1,\dots,A_N))$,
the latter being the measures on the finite $\sigma$-algebra generated by
$A_1,\dots,A_N$.

\subsection{Construction of the measure on the full algebra}

Note that generally a continuous function is not measurable with respect to 
the $\sigma$-algebra spanned by a finite partition $A_1,\dots,A_n$.
However, it can be arbitrarily well-approximated in the sup-norm
by such functions.
Specifically, because $\mathscr{X}$ is compact, we can find a triangular
family of sets $(A_{n,k}\, :\, n \in \N_+,\, 1\leq K\leq N(n))$,
where each $N(n) \in \N_+$ and such that each set $A_{n,k}$, for
$k=1,\dots,N(n)$,
has diameter no larger than $1/n$.
This is just by choosing finite subcovers of the cover of $\mathscr{X}$ by
all open sets of diameter less than $1/n$.
Any continuous function on a compact set is necessarily unifomly continuous
by the ``Heine-Cantor'' theorem.
Therefore, any continuous function can be arbitrarily well-approximated
by functions $f_n$ such that $f_n$ is constant on each
$A_{n,k}$, for $k=1,\dots,N(n)$
and $n \in \N_+$.
Thus, $f_n$ is measurable with
respect to the $\sigma$-algebra generated by
$(A_{m,k}\, :\, 1\leq m\leq n,\, 1\leq k\leq N(m))$.
Let us call the $\sigma$-algebra $\mathcal{G}_n$.

We can define a family of functions, $M^{(n)}$, such that 
\begin{equation}
  M^{(n)}\, :\, [0,1] \to \mathcal{M}_+^1(\mathcal{G}_n)
\end{equation}
as constructed in the last section,
choosing the Borel partition of $\mathscr{X}$ to correspond
to $\mathcal{G}_n$, 
(More specifically, 
we choose an order for the sets
$(A_{m,k}\, :\, 1\leq m\leq n,\, 1\leq k\leq N(m))$,
and take the difference of the $i$th set with the union of its
$i-1$ predecessors for each $i$.)
Almost surely, with respect to $\mu^{\infty}$,
and for almost every $t$, this function satisfies the conditions to
be a measure.
Moreover, it is trivial to see that the measures given by
$(M^{(n)}\, :\, n \in \N_+)$ are consistent, i.e., form a projective family
(for almost all $t$).
At this point, we could use Caratheodory's theorem to construct
$\mu(dx,t)$, but we will proceed differently.

Suppose that $f \in \mathscr{X}$ is any continuous function, and
$E \subset [0,1]$ is any Borel measurable set.
Note that each $f \in \mathcal{C}(\mathscr{X})$ is the uniform
limit of functions $f_n \in \mathcal{B}(\mathscr{X})$, (this space
is the set of everywhere-bounded, Borel measurable functions),
such that each $f_n$ is measurable with respect to the algebra
$\mathcal{G}_n$.
By H\"older's theorem, we can prove that 
\begin{equation}
  \int_{0}^1 \int_{\mathscr{X}} f_n(x)\, \mu_{M_n(t)}(dx)\, \mathbbm{1}_E(t)\, dt
\end{equation}
is a Cauchy sequence.
Let us define the limit as $\ell(\mathbbm{1}_E \otimes f)$.
Then it is easily seen that for each non-null $E$, the linear functional
$\ell(\mathbbm{1}_E \otimes \cdot)/|E|$, defines a measure,
$\nu_E \in \mathcal{M}^+(\mathscr{X})$,
by the Riesz-Markov theorem.
These measures are consistently defined, so as to define
a measure $\nu \in \mathcal{M}^+_1(\mathscr{X}\times [0,1])$
such that $\nu(E\times \mathscr{X}) = \nu_E(\mathscr{X}) = |E|$
for all Borel $E \subset [0,1]$.
Repeating the argument from Section 2, in which we constructed 
$\mu(dx,t)$ from
$\nu \in \mathcal{M}^+_1(\mathscr{X}\times [0,1])$
sastisfying $\nu(\mathscr{X}\times E)=|E|$ for all Borel $E\subset [0,1]$,
we can now construct the $\mu^{(\infty)}\otimes(\textrm{Lebesgue})$-almost
surely unique measure $\mu(dx,t)$ which satisfies
$\ell(\mathbbm{1}_E \otimes f)=\ip{\mu}{\mathbbm{1}_E\otimes f}$
for every $E\subset [0,1]$ and $f \in \mathscr{X}$.

We define $P(d\mu)$ to be the distribution of $\mu$, which is a 
probability measure, measurable with respect to the tail $\sigma$-algebra
$\mathcal{F}_\infty$.

\subsection{Reconstruction of $\mu_n$}
\label{reconstruct}

Let us consider, again, the case of a finite partition $(A_1,\dots,A_R)$.
Let $n \in \N_+$, and we ask for the probability
\begin{equation}
  \p(X^{(n)}_1 \in A_{i(1)},\dots,X^{(n)}_n \in A_{i(n)})\, ,
\end{equation}
for $i:\{1,\dots,n\} \to \{1,\dots,R\}$.
Using the projective property of the conditional expectation, the so-called ``law of the iterated
conditional expectation'', we obtain
\begin{equation}
\begin{split}
  \p(X^{(n)}_1 \in A_{i(1)},\dots,X^{(n)}_n \in A_{i(n)})\, \\
  &\hspace{-50pt} =\, \E\{\mathbbm{1}[X^{(n)}_1\in A_{i(1)},\dots,X^{(n)}_n \in A_{i(n)}]\} \\
  &\hspace{-50pt} =\, \E\{\phi_N[\mathbbm{1}[X^{(n)}_1\in A_{i(1)},\dots,X^{(n)}_n \in A_{i(n)}]]\}\, ,
\end{split}
\end{equation}
for any $N\geq n$.

We observe that
\begin{multline}
  \phi_N[\mathbbm{1}[X^{(n)}_1\in A_{i(1)},\dots,X^{(n)}_n \in A_{i(n)}]]\\
  =\, \sum_{1\leq j(1)<\dots<j(n)\leq N} \binom{N}{n}^{-1} \mathbbm{1}[X^{(N)}_{j(1)}\in A_{i(1)}] \cdots \mathbbm{1}[X^{(N)}_{j(n)} \in A_{i(n)}]\, .
\end{multline}
Defining the density $p(N,n;\cdot) : \Delta^{n-1} \to \R_+$ by 
\begin{equation}
\begin{split}
  &p(N;n;t_1,\dots,t_n)\\ 
  &=\,  \sum_{1\leq j_1<\dots<j_n\leq N} \binom{N}{n}^{-1} \mathbbm{1}\left[\frac{j_1-1}{N}\leq t_1<\frac{j_1}{N}\right] \cdots 
  \mathbbm{1}\left[\frac{j_n-1}{N}\leq t< \frac{j_n}{N}\right]\, ,
\end{split}
\end{equation}
we see that, in the vague topology,
\begin{equation}
  \lim_{N \to \infty} p(N,n;t_1,\dots,t_n)\, dt_1\cdots dt_n\, =\, n!\, dt_1\cdots dt_n\, ,
\end{equation}
 on $\Delta^{n-1}$, and also that 
\begin{equation}
\begin{split}
&\phi_N[\mathbbm{1}[X^{(n)}_1\in A_{i(1)},\dots,X^{(n)}_n \in A_{i(n)}]]\\
&=\, \int_{\Delta^{n-1}} \left[\prod_{k=1}^n \int_{\mathscr{X}} \mathbbm{1}[x_k \in A_{i(k)}]\, M_N(dx_k,t_k)\right]\, 
   p(N,n;t_1,\dots,t_n)\, dt_1\cdots dt_n\, .
\end{split}
\end{equation}
Note that this last function is obviously continuous in the weak $\mathcal{L}$-topology,
since it can be approximated in norm by a sum of products of linear functions in $\mathcal{L}$.

Then by the vague convergence of the density, as well as the convergence of $M_n$ in the weak
$L^1$-topology, we conclude that, $\mu^{(\infty)}$-almost surely
\begin{multline}
  \phi_\infty[\mathbbm{1}[X^{(n)}_1\in A_{i(1)},\dots,X^{(n)}_n \in A_{i(n)}]]\\
  \qquad =\, \int_{\Delta^{n-1}}  \left[\prod_{k=1}^n \int_{\mathscr{X}} \mathbbm{1}[x_k \in A_{i(k)}]\, M(dx_k,t_k)\right]\, 
   n!\, dt_1\cdots dt_n\, .
\end{multline}
This is exactly what we needed to prove in order to demonstrate that
\begin{multline}
  \p(X^{(n)}_1 \in A_{i(1)},\dots,X^{(n)}_n \in A_{i(n)})\\
  =\, \int_{\mathcal{M}^+_1(\mathscr{X})} \int_{\Delta^{n-1}} \left[\prod_{k=1}^n \int_{\mathscr{X}} \mathbbm{1}[x_k \in A_{i(k)}]\, \mu(dx_k,t_k)\right]\,
  n!\, dt_1\cdots dt_n\, P(d\mu)\, .
\end{multline}
The result for continuous functions follows by approximation.

\subsection{Uniqueness of the measure}
The uniqueness of the measure will follow by a standard weak LLN argument.
The proof of the weak LLN is also standard, obtained by bounding
second moments.
However, to apply weak LLN, we need one more topological fact,
which is that the algebra of continuous functions
on $\mathcal{K}$, spanned by products of functions
in $\mathcal{L}$, is dense in $\mathcal{C}(\mathcal{K})$.
Since $\mathcal{K}$ is a compact metric space (as noted in Section 2), 
this fact will
follow by the Stone-Weierstrass theorem, if we verify
that $\mathcal{L}$ separates points in $\mathcal{K}$.

There are various ways to prove that $\mathcal{L}$ separates
points in $\mathcal{K}$.
We choose the following.

Let $\mathcal{M}(\mathscr{X})$ be the Banach space
of all signed measures on $\mathscr{X}$ with finite
total-variation norm 
\begin{equation}
  \|\mu\|_{TV}\, =\, \int_{\mathscr{X}} |\mu|(dx)\, .
\end{equation}
Define $\widetilde{\mathcal{K}} = L^\infty([0,1],\mathcal{M}(\mathscr{X}))$
with the norm
\begin{equation}
  \|\mu\|\, =\, \textrm{ess}\sup \{ \|\mu(\cdot,t)\|_{TV}\, :\, t \in [0,1]\}\, .
\end{equation}
This space is not separable, and generally the topology is 
quite different from that of $\mathcal{K}$.
(More directly, this space is horrible from the perspective of
doing analysis on it.)
However, note that $\mathcal{K}$ is 
definitely a convex subset of $\widetilde{\mathcal{K}}$.

For any $\mu \in \widetilde{\mathcal{K}}$, we observe that
\begin{equation}
  F \mapsto \ip{\mu}{F}\, =\, \int_0^1 \int_{\mathscr{X}} F(x,t)\, \mu(dx,t)\, dt
\end{equation}
defines a continuous linear functional on $\mathcal{L}$,
in fact $|\ip{\mu}{F}| \leq \|\mu\|\cdot \|F\|$, by H\"older's
inequality.
Also, by the usual argument, (c.f.\ \cite{LiebLoss}, Theorem 2.10),
the set of linear functionals, $\{\ip{\mu}{\cdot}\, :\, \mu \in \widetilde{\mathcal{K}}\}$,
is separating
in $\mathcal{L}$.
I.e., for every nonzero $F \in \mathcal{L}$, there is a $\mu \in \widetilde{\mathcal{K}}$
such that $\ip{\mu}{F} \neq 0$.
Therefore, $(\mathcal{L},\widetilde{\mathcal{K}})$
forms a ``dual pair'' (c.f.\ \cite{Simon}, Section 1.5).

We observe that $\mathcal{K}$ is a convex, closed subset
of $\widetilde{\mathcal{K}}$ in the weak-$\mathcal{L}$
topology.
Therefore, by an application of 
the separating hyperplane theorem, one can prove that
$\mathcal{L}$ does separate points in $\mathcal{K}$,
(c.f.\ \cite{Simon}, Theorem 1.5.2, which states
that for any relatively closed subset $K \subset \mathcal{K}$
and any point $\mu \in \mathcal{K}\setminus K$,
there is a linear functional equal to 1 on $K$ and strictly
less than 1 at $\mu$).

Since $\mathcal{L}$ separates points by itself, it is clear that
the algebra generated by $\mathcal{L}$ satisfies the hypotheses
of the Stone-Weierstrass theorem.
As a final detail, we observe that taking $F(x,t) \equiv 1$
yields the constant function 1 on $\mathcal{K}$.
Therefore, using the Riesz-Markov theorem and the fact
that $\mathcal{K}$ is a compact metric space,
we can reduce the task of proving that two probability measures
$P_1, P_2 \in \mathcal{M}^+_1(\mathcal{K})$
are identical, to the easier task of proving that for every function
in the algebra generated by $\mathcal{L}$, we obtain
the same expectation.
This is what we show next by using the weak law of large numbers.

Let $f_1,\dots,f_n \in \mathcal{C}(\mathscr{X})$
and let $I_1,\dots,I_R \subset [0,1]$ be disjoint open intervals.
We intend to show that
\begin{multline}
\label{wlln}
  \lim_{n \to \infty}
  \int_{\mathscr{X}^n} \mu_n(dx_1\otimes\cdots\otimes dx_n)\, 
  \prod_{r=1}^R \left(\sum_{k = 1}^n n^{-1}\, \mathbbm{1}[(k/n) \in I_r] f_r(x_k) \right)\\
  =\, \int_{\mathcal{M}^+_1(\mathscr{X})} P(d\mu)\,
  \prod_{r=1}^R \left(\int_{I_r} \int_{\mathscr{X}} f_r(x_k)\, \mu(dx_r,t_r)\, dt_r\right)\, .
\end{multline}
Because of the remarks above about the Stone-Weierstrass theorem,
and because Lebesgue measure is regular, this will prove
that $P$ is uniquely determined by $(\mu_n\, :\, n \in \N_+)$,
by an approximation argument.

We will omit some of the details needed to prove (\ref{wlln}).
It will suffice to consider the one- and two-coordinate marginals 
of the Lebesgue measure on $\Delta^{n-1}$.
Specifically, that
\begin{equation}
\begin{split}
  &n!\, \int_{[0,1]^{k-1}} \mathbbm{1}[0<t_1<\dots<t_{k-1}<t]\, dt_1\cdots dt_{k-1}\\
  &\times \int_{[0,1]^{n-k}} \mathbbm{1}[t<t_{k+1}<\dots<t_n<1]\,
  dt_{k+1}\cdots dt_n\\
  &\hspace{150pt} =\, \frac{t^{k-1} (1-t)^{n-k}\, \mathbbm{1}[0<t<1]}{B(k,n-k+1)}\, ,
\end{split}
\end{equation}
and that 
\begin{equation}
\begin{split}
  &n!\, \int_{[0,1]^{j-1}} \mathbbm{1}[0<t_1<\dots<t_{j-1}<s\cdot t]\, dt_1\cdots dt_{j-1}\\
  \qquad &\times \int_{[0,1]^{k-j-1}} \mathbbm{1}[s\cdot t<t_{j+1}<\dots<t_{k-1}<t]\,
  dt_{j+1}\cdots dt_{k-1}\\
  \qquad &\times \int_{[0,1]^{n-k}} \mathbbm{1}[t<t_{k+1}<\dots<t_{n}<1]\,
  dt_{k+1}\cdots dt_{n}\\
  &=\, \frac{s^{j-1} (1-s)^{k-j-1}\, \mathbbm{1}[0<s<1]}{B(j,k-j)}\cdot
  \frac{t^{k-1} (1-t)^{n-k}\,  \mathbbm{1}[0<t<1]}{B(k,n-k+1)}\, .
\end{split}
\end{equation}
Using this, we can show that when $(t_1,\dots,t_n)$ are distributed uniformly on $\Delta^{n-1}$,
we have
\begin{equation}
  \E\{t_k\}\, =\, \frac{k}{n+1}\qquad \mbox{and, for $j\leq k$,}\qquad
  \E\{t_j t_k\}\, =\, \frac{j(k+1)}{n(n+1)}\, .
\end{equation}
This is enough to show that the empirical measures
\begin{equation}
  \nu_n^{(r)}(\cdot;t_1,\dots,t_n)\, =\, n^{-1} \sum_{k=1}^n \mathbbm{1}[(k/n) \in I_r] \delta_{t_k}\, ,
\end{equation}
are asymptotically independent for $r\neq s$, as $n \to \infty$,
and that in the vague topology
\begin{equation}
\label{lim}
  \lim_{n \to \infty} \nu_n^{(r)}(dt;t_1,\dots,t_n)\, =\, \mathbbm{1}[t \in I_r]\, dt\, ,
\end{equation}
almost surely.

Using this, one can see that
\begin{multline}
  \lim_{n \to \infty} \E\left\{\prod_{r=1}^R \left(\sum_{k = 1}^n n^{-1}\, \mathbbm{1}[(k/n) \in I_r] f_r(X^{(n)}_k) \right)\right\}\\
  =\, \int_{\mathcal{M}^+_1(\mathscr{X})}P(d\mu)\, \prod_{r=1}^R \left(\int_{I_r} \int_{\mathscr{X}} f_r(x_r)\, \mu(dx_r,t_r)\, dt_r\right)\, .
\end{multline}
For example, for $n \in \N_+$, one can consider the reversed
martingale based on the $n$th term.
Using calculations as in subsection \ref{reconstruct}, one can see that
the $\mu^{\infty}$-almost sure limit of the reversed semimartingale based on the $n$th term is
\begin{equation}
  \int_{\Delta^{n-1}} \prod_{r=1}^R \left(\int_0^1 \int_{\mathscr{X}} f_r(x_k)\, \mu(dx_r,s_r)\, \nu_n^{(r)}(ds_r;t_1,\dots,t_n)\right) 
  n!\, dt_1\cdots dt_n\, .
\end{equation}
Indeed, this relies on nothing other than the restriction property of the Lebesgue measure
on the simplex:
for $N\geq n$,
\begin{multline}
  \sum_{1\leq k_1<\dots<k_n\leq N} \binom{N}{n}^{-1} \int_{\Delta^{N-1}} \delta(s_1-t_{k_1})\cdots \delta(s_n-t_{k_n})\,
  N! dt_1\cdots dt_N\\ 
  =\, n! \mathbbm{1}[0<s_1<\dots<s_n<1]\, .
\end{multline}
Note that this property is obvious in the probabilistic framework, where the law of $(t_1,\dots,t_N)$ is obtained
by randomly choosing $N$ points, uniformly and independently in $[0,1]$, and then sorting.

Using the asymptotic independence of the empirical measures, and the limit of (\ref{lim}),
one deduces equation (\ref{wlln}).
This completes the proof of the theorem.

\section{Classical Spin Systems}
\label{Sec:CSS}

Our first application is to classical spin chains with
a slightly unusual scaling of the interaction.
We will explain why one might be interested in these models, later.

For simplicity, we consider the case of two-valued spins.
I.e., $\mathscr{X} = \{+1,-1\}$.
We will consider a Hamiltonian on $N$ spins, with a fixed $n$-body interaction.
Afterwards, we will consider the thermodynamic limit $N \to \infty$.
The $n$-body interaction is supposed to be defined by a function
\begin{equation}
  \phi_n : \mathscr{X}^n \to \R\, .
\end{equation}
We will be most interested in the case that $\phi_n$ is an asymmetric function
of its variables; i.e., that the order of the variables
$(\sigma_1,\dots,\sigma_n)$ does matter.

For each $N \in \N_+$, we let $\Lambda_N = \{1,\dots,N\}$
be the one-dimensional spin system.
For $k\in\{0,\dots,N\}$, let $\mathcal{P}_k(\Lambda_N)$
be the collection of all cardinality-$k$ subsets of $\Lambda_N$.
For any $\sigma \in \mathscr{X}^N$ and any $X \subset \Lambda_N$,
let $\sigma \restriction X$ be the restrition of $\sigma$ to $X$,
so that $(\sigma \restriction X) \in \mathscr{X}^X$.
For $N\geq n$, define the Hamiltonian $H_N : \mathscr{X}^N \to \R$ as
\begin{equation}
  H_N(\sigma)\, =\, N\, \sum_{X \in \mathcal{P}_n(\Lambda_N)}
  \binom{N}{n}^{-1}\, \phi_n(\sigma \restriction X)\, .
\end{equation}
For any $X \in \mathcal{P}_n(\Lambda_N)$, we keep the relative order
of the coordinates.
We would like to call these models ``oriented, mean-field models''.

Let us now discuss the context of the present models, within the 
framework of short-ranged models, and exchangeable, mean-field models.
Exchangeable, mean-field models of classical statistical mechanics,
such as the Curie-Weiss model, have long been considered as 
simplified analogues of more realistic short-ranged models.
However, while there were heuristic links between mean-field models
and short-ranged model,
they remained logically disconnected until
the work of Lebowitz and Penrose \cite{LP}.
(Also see the reference\footnote{
Lebowitz and Penrose considered a continuum mean-field model,
that of the van der Waals gas. In \cite{Thompson}, Thompson gives
the application of the Lebowitz-Penrose theorem to the Curie-Weiss
model, which is simpler.
So, for pedagogical reasons, one might look at that first.}
\cite{Thompson}, Appendix C.)

Lebowitz and Penrose showed how to recover mean-field models,
along with the proper Maxwell construction for the order-parameter
in the presence of a phase transition.
Specifically, they considering a range-$R$ model on a lattice
of linear size $N$. 
Lebowitz and Penrose took $R\to \infty$
after taking $N \to \infty$.
In the present section, we are considering another asymptotic
regime, where $R$ and $N$ go to $\infty$, together.
In this case, the exchangeable scaling of Lebowitz and Penrose
may be seen as an infinitesimal $t$ interval,
and if $P$ is concentrated on $\mu(dx,t)$ which
are continuous with respect to $t$, 
then we could say that those 
thinning-invariant states are ``locally exchangeable''.

We are primarily interested in the Gibbs measures.
Suppose that $\beta \in [0,\infty)$.
Then the Gibbs state for $\beta H_N$ is 
\begin{equation}
  \rho_{N,\beta}(\{\sigma\})\, =\, Z_N(\beta)^{-1}\, 2^{-N}\, \exp(-\beta H_N(\sigma))\, ,
\end{equation}
where $Z_N(\beta)$ is the normalized partition function
\begin{equation}
  Z_N(\beta)\, =\, \sum_{\sigma \in \mathscr{X}^N}\, 2^{-N}\, \exp(-\beta H_N(\sigma))\, .
\end{equation}

Given $\rho_{N,\beta}$, which is a measure on $\mathscr{X}^N$, we can define a measure
$\rho^{(N)}_k \in \mathcal{M}^+_1(\mathscr{X}^k)$ for all $k\leq N$, by
\begin{equation}
  \rho^{(N)}_k(\{\sigma^{(k)}\})\, 
  =\, \sum_{X \subset \mathcal{P}_k(\Lambda_N)} \binom{N}{k}^{-1}
  \sum_{\sigma \in \mathscr{X}^N} \mathbbm{1}[\sigma \restriction X\, =\, \sigma^{(k)}]\, 
  \rho_{N,\beta}(\{\sigma\})\, .
\end{equation}
This is defined in a thinning-invariant way for $k\leq N$.
It is evident that every limit point of such sequences, as $N \to \infty$,
is also thinning-invariant.
By the Banach-Alaoglu theorem, there is always at least one convergent subsequence.
By Theorem 1, we have a representation for all such limit
points.

Of course, we could do this procedure for any sequence of finite spin chains,
not just the oriented, mean-field models we have defined here.
But the limit of such thinning-invariant states
will be irrelevant for short-ranged models, whereas
in the present situation, the limit points are relevant.
This is because $H_N$ is measurable with respect to $\rho^{(N)}_n$.
In fact,
\begin{equation}
  N^{-1}\, \rho_{N,\beta}(H_N)\, =\, n^{-1}\, \rho^{(N)}_n(\phi_n)\, ,
\end{equation}
where by this notation we mean to take expectations.

We would like to calculate the limit points of the Gibbs distribution
using our main theorem, instead of merely deducing a rather general
representation for them.
At the very least, one might expect to be able to prove that the extreme points
of the limit Gibbs distributions are also extreme points in the simplex of 
thinning-invariant sequences of measures.  
At the present time this is beyond our reach for general interactions, 
in part because we do not know
any intrinsic definition of limit Gibbs states, akin to the 
Dobrushin-Lanford-Ruelle condition for short-ranged models.
(It should be noted that the same problem exists for 
exchangeable, mean-field models like the Curie-Weiss model.)
We discuss this a bit more at the end of Section 7.

Since $H_N$ is measurable with respect to 
$\rho^{(N)}_n$,
we can at least determine the ground states and ground state energy density
in this way.

\subsection{Example 1: Curie-Weiss model with inhomogeneous field}

Let us take the case of a two-body interaction
\begin{equation}
  - \phi_2(\sigma_1,\sigma_2)\, =\, \frac{1}{2} \sigma_1 \sigma_2 + b \cdot (\sigma_2 - \sigma_1)\, ,
\end{equation}
By reflecting $\Lambda_N$ about $(N+1)/2$, we can flip $b$ for $-b$, so let us assume
that $0\leq b<\infty$.
Then, for $N\geq 2$,
\begin{equation}
  - H_N(\sigma)\, 
  =\, \frac{1}{N-1} \sum_{1\leq j<k\leq N} \sigma_j \sigma_k\, 
  +\, b \sum_{j=1}^N \frac{2j-N-1}{N-1}\, \sigma_j\, .
\end{equation}
If we restrict attention to extreme points of the thinning-invariant measures,
which we know must optimize the {\em linear} functional $H_N$, then 
we see that
\begin{equation}
  - \mu_2(\phi_2)\, 
  =\, \frac{1}{2} \left(\int_0^1 \varphi(t)\, dt\right)^2\, 
  +\, 2 b \int_0^1 (2t-1)\, \varphi(t)\, dt\, 
\end{equation}
where we define $\varphi(t) = \E^{\mu(\cdot,t)}\{\sigma\} = \mu(\{+1\},t) - \mu(\{-1\},t) \in [-1,1]$.

Let us condition on the value of 
\begin{equation}
\overline{m}\, =\, \int_0^1 \varphi(t)\, dt\, .
\end{equation}
Then we have to minimize the value of the linear functional
\begin{equation}
  \int_0^1 (2t-1)\, \varphi(t)\, dt\, 
\end{equation}
subject to the constraint on $\overline{m}$,
and, of course, $-1\leq \varphi\leq 1$.
One can easily conclude that the constrained optimizer is
\begin{equation}
  \varphi(t)\, =\, \begin{cases} -1  & \mbox{for } 0\leq t<t_0\, ,\\ 
    +1 & \mbox{for } t_0\leq t<1\, ,
\end{cases}
\end{equation}
where $t_0 = (1-\overline{m})/2$.
For example, one can consult \cite{LiebLoss}, Theorem 1.14, for a proof of 
the so-called ``bathtub principle''.

So the optimum value of this term, subject to the constraint $\overline{m}$,
is 
\begin{equation}
  \int_0^1 (2t-1)\, \varphi(t)\, dt\, 
  =\, 2t_0(1-t_0)\, 
  =\, (1-\overline{m}^2)/2\, .
\end{equation}
Therefore, the total energy is 
\begin{equation}
  -\mu_2(\phi_2)\, =\, \frac{b}{2}\, +\, \frac{1-b}{2} \overline{m}^2\, .
\end{equation}

In this case, we see that for $b<1$, there are two distinct solutions:
$\mu(\cdot,t) \equiv \delta_{+1}$ and $\mu(\cdot,t) \equiv \delta_{-1}$.
If $b>1$, there is a unique solution, such that
\begin{equation}
  \mu(\cdot,t)\, =\, \begin{cases}\delta_{-1} & 0<t<1/2\, ,\\
      \delta_{+1} & 1/2<t<1\, .
\end{cases}
\end{equation}
The most interesting situation is $b=1$.
Then the interaction is, modulo
a constant shift, $\phi_2(\sigma_1,\sigma_2) = (\sigma_1+1)(\sigma_2-1)$.
This is the mean-field version of the Ising ``kink'' Hamiltonian.
In the thermodynamic limit, a continuum of ground states exist:
for any $s \in [0,1]$, the measure
\begin{equation}
  \mu_s(\cdot,t)\, =\, \begin{cases}\delta_{-1} & 0<t<s\, ,\\
      \delta_{+1} & s<t<1\, .
\end{cases}
\end{equation}
is a ground state.
This is the continuum limit of the Ising kink ground states.

It may be noted that for this simple model, one can calculate 
the ground-states explicitly, for every $n \in \N_+$,
thus confirming our conclusions.

But even for more complicated models, where one cannot calculate
the finite volume ground states, we know rigorously that any 
limit point must satisfy the conclusions of our theorem.
Therefore, we can calculate the limiting ground state energy density,
as well as all possible limit points of ground states,
using the same technique as above (assuming that we can solve the
continuous optimization problem).

\section{Exclusion processes}

Another interesting application of Theorem 1 is for limits
of invariant measures for exclusions processes, with the same
scaling of the interaction range as in Section 4.

For these models, $\mathscr{X} = \{0,1\}$.
Suppose the set of sites, $S$, is finite.
A particle configuration is denoted $\eta:S\to \mathscr{X}$.
Then, following reference \cite{Liggett}, Chapter VIII,
the simple exclusion process is defined to be the Feller
semigroup with generator given by $\Omega$:
\begin{equation}
  \Omega f(\eta)\, =\, \sum_{x,y \in S}\, \mathbbm{1}[\eta(x)=1,\, \eta(y)=0]\,
  \left(f(\eta_{xy}) - f(\eta)\right)\, p(x,y)\, ,
\end{equation}
where for each pair $x,y \in S$, $p(x,y)\geq 0$ (the numbers $p(x,x)$ are 
obviously irrelevant), and for any $x\neq y$, the configuration
$\eta_{xy}$ is obtained from $\eta$ by switching the values
at $x$ and $y$:
\begin{equation}
  \eta_{xy}(u)\, =\, \begin{cases} \eta(y) & \mbox{if } u=x\, ,\\
    \eta(x) & \mbox{if } u=y\, ,\\
    \eta(u) & \mbox{otherwise.} \end{cases}
\end{equation}

It is clear that $\Omega$ is a positivity-preserving form,
therefore by the Perron-Frobenius theorem, there is a unique
positive eigenvector in each ergodic sector of $\Omega$.
The ergodic sectors (or components) are the maximal sets of 
configurations which are 
connected in the graph obtained by connecting every pair of 
configurations corresponding to a nonzero matrix entry.
Assuming that for every $x\neq y$ there is a $k \in \N_+$,
and a sequence $x=x_1,\dots,x_k=y$ with $p(x_j,x_{j+1})>0$
for $j=1,\dots,k-1$, it will be clear that the connected
components are the sets of configurations with a fixed
total number of particles.
(There would be subtleties if a configuration $\eta$
could be connected to $\zeta$, but not $\zeta$ to 
$\eta$, and this eliminates that.)
I.e., the generator of the dynamics commutes with the total
particle number, but every pair of configurations
with the same number of particles are eventually connected
by $\Omega$.

We will consider the case that, for each $N \in \N_+$,
the sites are $S_N  = \{1,\dots,N\}$ and 
\begin{equation}
  p_N(x,y)\, =\, \begin{cases} q/(N-1) & \mbox{if } x<y\, ,\\ (1-q)/(N-1) & \mbox{if } y<x\, . \end{cases}
\end{equation}
We will consider the case $1/2<q<1$, the opposite case being treated
symmetrically.
Note that in this case $\Omega_N$ is not symmetric.
Also, for $N>2$, the matrix $p_N(x,y)$ has no reversible measure.
I.e., there is no $\pi : S_N \to \R_+$ such that
$\pi(x) p_N(x,y) = \pi(y) p_N(y,x)$ for all $x\neq y$.
That would require $\pi(x)/\pi(y) = (1-q)/q$ whenever
$x<y$ which is clearly impossible when $1/2<q<1$. 
The reason this is important is that Liggett and others
have made a rather careful analysis of the complete set of invariant
measures for the exclusion process when such a $\pi$ exists.
(See \cite{Liggett}, Chapter VIII; \cite{Liggett2}, Part III;
and \cite{Jung}.)
To the best of our knowledge, the present sequence of models has not been 
solved, previously.

Letting $\mu_{N,n}$ be the invariant measure for the size $N$ system
with $n$ particles, we can again construct $\mu^{(N,n)}_k$
for $1\leq k\leq N$ obtained by choosing $k$ sites at random
and projecting the measure to cylinder sets based at those sites,
then averaging over the $k$ sites chosen.
It should be noted that the measure $\mu^{(N,n)}_k$ does not
have a pure value of total number of particles,
but on average the density is $n/N$.
If we take the limit $n,N \to \infty$,  such that 
$n/N \to \rho \in [0,1]$, we know that there will be a weak-$*$
limit point of the thinning-invariant triangular
array, and also that any limit point will be thinning invariant.

As usual, we do not want to merely know that a representation of
such limit points exists, we actually want to calculate the limit
point.
In the example which follows we will {\it assume} that there is a 
limit point which is extremal in the simplex of thinning-invariant
measures.
We would like to point out that this is a non-void assumption:
we know of no inherent reason that the limit point should be contained
in the set of extreme points of the thinning-invariant simplex
\footnote{We do not rule out the possibility that a theorem which we are unaware
of does guarantee that the limit will be extremal.
For example, for translationally-invariant, short-ranged classical spin systems, 
the extreme translation-invariant limit Gibbs states
are extreme in the family of all translation-invariant states, because the limiting
entropy density is affine on that simplex, even though its finite-dimensional approximations are not.
C.f., \cite{Israel} or \cite{Simon}.}.
For, although $\Omega_N$ is a linear operator,
one does not determine the invariant measure by optimizing the
expectation of $\Omega_N$.
Rather, one optimizes the lower Collatz-Wielandt number 
\begin{equation}
  l(\Omega_N,\mu)\, =\, \max\left\{\lambda \in \R\, :\, ([\Omega_N]_* - \lambda)\mu \geq^{K} 0\right\}\, ,
\end{equation}
where $\mu$ is restricted to be a nonnegative measure and $\geq^K$ is the order
on $\mathcal{M}(\mathscr{X}^N)$ 
with respect to the cone of nonegative measures $K = \mathcal{M}^+_1(\mathscr{X}^N)$.

This is a nonlinear function of $\mu$, therefore, there is no reason that it will be attained
on an extremal (i.e., delta-) measure.

We will, nonetheless, find a condition on extreme points of the thinning-invariant
simplex $P = \delta_{\mu}$, which is necessary to be a limit of invariant measures
for the simple exclusion process.
The condition, at the very least, leads to a heuristic derivation of a limiting
invariant measure;
however it should be noted that the heuristic assumption, that the limiting measure
factorizes, is completely unoriginal.
For example, it is this assumption, along with some elementary analysis,
which leads to the Burger's equation for the limiting invariant measures
of the short-ranged ASEP in \cite{Liggett2}, pp.\ 222-224.

What is new is to connect this heuristic assumption with a rigorous
theorem: that if the limit of the invariant measures would be an extreme
point of the thinning-invariant simplex, then it {\em must} satisfy
the following conditions.
As we will see, the conditions are sufficient to uniquely specify one
extreme point for each choice of $\rho \in [0,1]$.

\subsection{Example 2: Mean-field ASEP}

Note that, if we define the evaluation maps $\epsilon_x(\eta) = \eta_x$, then
we may rewrite
\begin{equation}
  \Omega_N f(\eta)\, =\, \sum_{x,y \in S}\, \epsilon_x(\eta)\, [1-\epsilon_y(\eta)]\, 
  \left(f(\eta_{xy}) - f(\eta)\right)\, p_N(x,y)\, .
\end{equation}
If we consider $f = \epsilon_x$, then we find
\begin{equation}
\label{firsteq}
\begin{split}
  -\Omega_N \epsilon_x\, &=\, 
   \sum_{y<x} \left( \frac{1-q}{N-1}\, \epsilon_x\, [1-\epsilon_y]\, -\, \frac{q}{N-1}\, \epsilon_y\, [1-\epsilon_x]\right)\\
  &\qquad + \sum_{y>x} \left( \frac{q}{N-1}\, \epsilon_x\, [1-\epsilon_y]\, -\, \frac{1-q}{N-1}\, \epsilon_y\, [1-\epsilon_x]\right)\, .
\end{split}
\end{equation}

Given $L \in \N_+$ and $\ell=1,\dots,L$, let us define
\begin{equation}
  f^{(N)}_{L,\ell}\, =\, \sum_{x \in S_N} \frac{1}{N}\cdot \mathbbm{1}\hspace{-2pt}\left[\frac{\ell-1}{L} < \frac{x}{N} \leq \frac{\ell}{L}\right]\, \epsilon_x\, .
\end{equation}
If $\mu_{N,n}$ does converge to a thinning-invariant measure which is extreme, i.e.,
given by a measure $\mu \in \mathcal{K}$ instead of a mixture $P \in \mathcal{M}^+_1(\mathcal{K})$,
then asymptotically $(f^{(N)}_{L,\ell}\, :\, \ell=1,\dots,L)$ forms an independent
family and
\begin{equation}
  \lim_{N \to \infty} \E^{\mu_N}\left\{\varphi(f^{(N)}_{L,\ell})\right\}\, =\, \varphi\left(\int_0^1 
  \mathbbm{1}\hspace{-2pt}\left[\frac{\ell-1}{L} < t \leq \frac{\ell}{L}\right]\, \mu(\{1\},t)\, dt\right)\, ,
\end{equation}
for any continuous function $\varphi$, by the weak LLN which we proved in Section 3.6.
Let us define $\phi(t) = \mu(\{1\},t)$ which exists as a function in $L^1(dt)$.

On the other hand, from (\ref{firsteq}), we can easily determine that 
\begin{equation}
\begin{split}
  -\Omega_N f^{(N)}_{L,\ell}\, 
  &=\, \sum_{\ell'<\ell} \left(\frac{1-q}{L-1}\, f^{(N)}_{L,\ell}\, \left[1-f^{(N)}_{L,\ell'}\right]\, -\, 
  \frac{q}{L-1}\, f^{(N)}_{L,\ell'}\, \left[1-f^{(N)}_{L,\ell}\right]\right)\\
  &\qquad +\,  \sum_{\ell'>\ell} \left(\frac{q}{L-1}\, f^{(N)}_{L,\ell}\, \left[1-f^{(N)}_{L,\ell'}\right]\, -\, 
  \frac{1-q}{L-1}\, f^{(N)}_{L,\ell'}\, \left[1-f^{(N)}_{L,\ell}\right]\right)\\
  &\qquad + O(1/L)\, .
\end{split}
\end{equation}
In an invariant measure, the expectation of the right-hand-side must equal zero, because $[\Omega_N]_* \mu_N = 0$.

Taking the ordered limit, where first $N \to \infty$, and then $L \to 0$, and using the factorization
property of the limit, we determine that a necessary requirement for the limit
to be described by an extreme thinning-invariant measure parametrized by $\mu$ is that
\begin{equation}
\begin{split}
\label{inteq}
  0\, &=\, \int_0^t \left((1-q)\, \phi(t)\, [1-\phi(s)]\, -\, q\, \phi(s)\, [1-\phi(t)]\right)\, ds\\
  &\qquad + \int_t^1 \left(q\, \phi(t)\, [1-\phi(s)]\, -\, (1-q)\, \phi(s)\, [1-\phi(t)]\right)\, ds\, ,
\end{split}
\end{equation}
for almost every $t \in [0,1]$.

Since the particle number is fixed by $\Omega_N$, we can fix the density
\begin{equation}
  [0,1] \ni \rho\, =\, \int_0^1 \phi(t)\, dt\, .
\end{equation}
We also define the function
\begin{equation}
  u(t)\, =\, \frac{1}{2} \left(\int_0^t \phi(s)\, ds\, -\, \int_t^1 \phi(s)\, ds\right)\, .
\end{equation}
Evidently, this is an absolutely continuous function, such that
\begin{equation}
\label{phitopsi}
  u(0)\, =\, -\rho/2\, ,\quad
  u(1)\, =\, \rho/2\, ,\quad \mbox{and}\quad
  u'\, =\, \phi\, ,\ \mbox{a.e.}
\end{equation}

\begin{figure}
\begin{center}
\resizebox{12truecm}{!}{\includegraphics{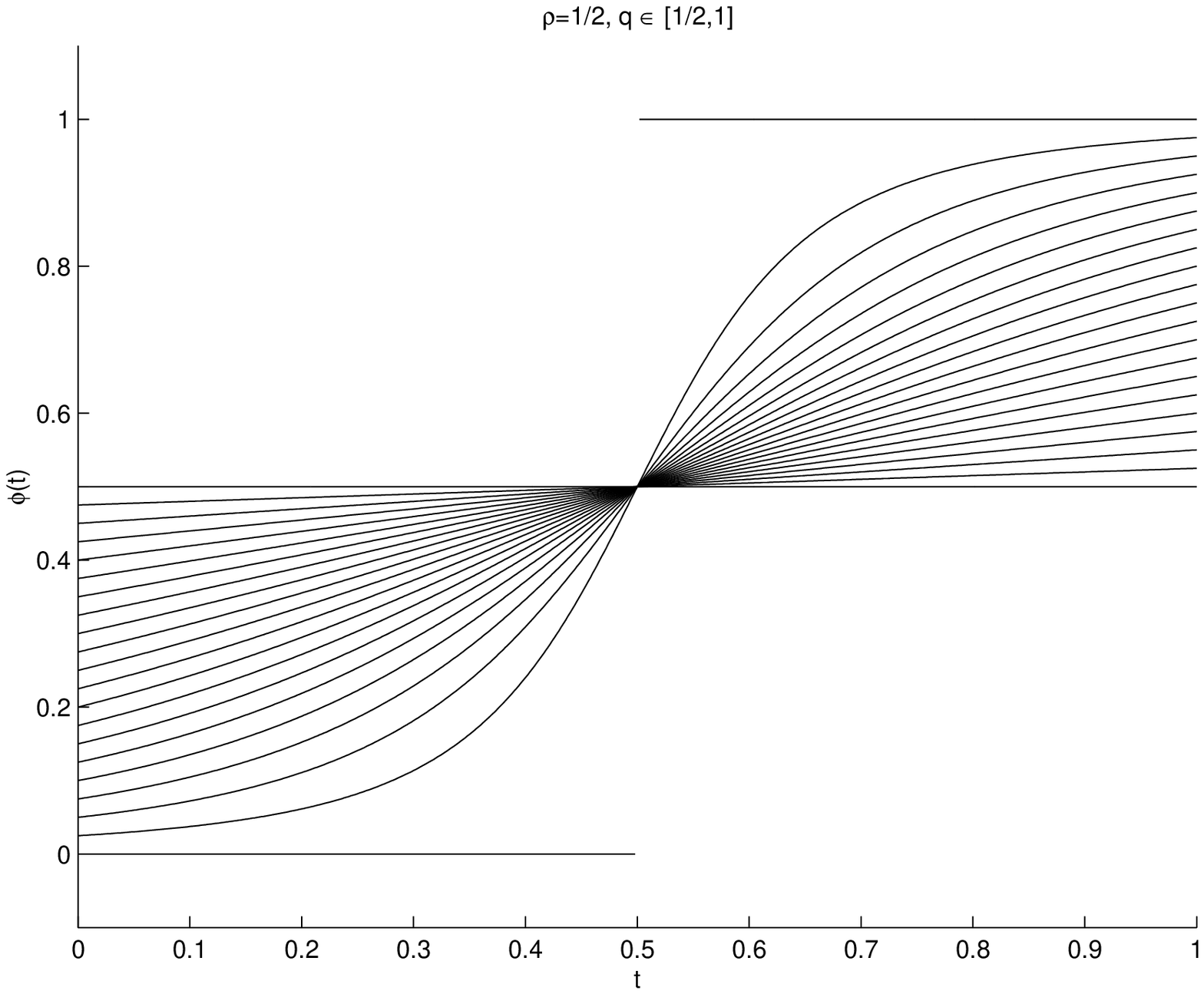}} \\
\resizebox{12truecm}{!}{\includegraphics{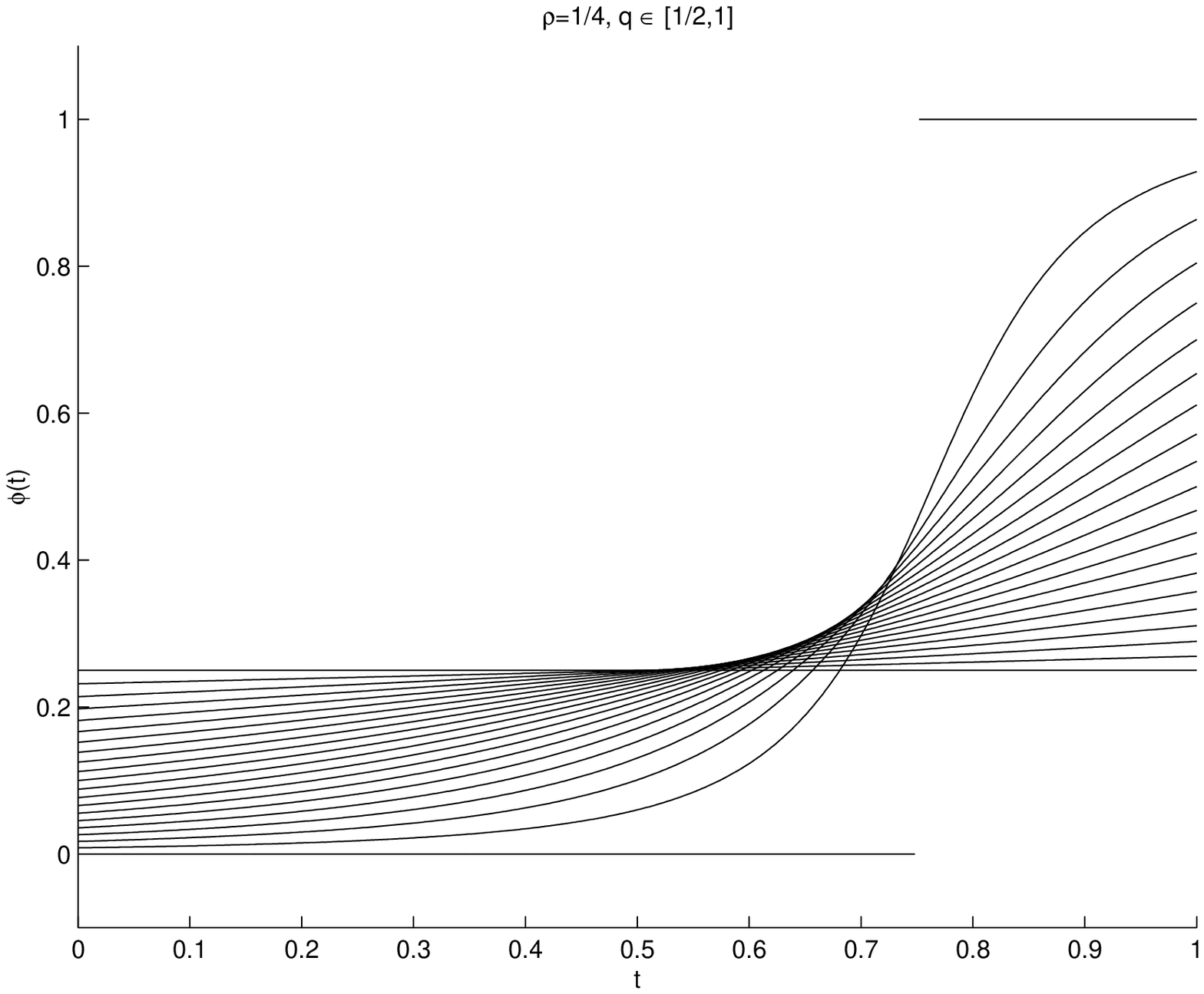}} \\
\end{center}
\caption{Density profile for an extreme, thinning-invariant
limit point
of invariant measures for the mean-field ASEP.}
\end{figure}

By elementary algebra, the integral equation (\ref{inteq}) 
can be rewritten as
\begin{equation}
\label{secdeq}
  2 \phi(t)\, =\, \frac{\frac{1}{2}\rho\, +\, (2q-1) u(t)}{\frac{1}{2}[(1-q)t+q(1-t)] + (2q-1)u(t)}\, .
\end{equation}
This shows that $\phi$ has an absolutely continuous version,
because $u$ is absolutely continuous, and the fraction does not diverge.
Using (\ref{phitopsi}), one can rewrite the almost everywhere identity (\ref{inteq})
in the even better form
\begin{equation}
  0\, =\, \frac{d}{dt} \left((2q-1) u(t)^2\, +\, [(1-t)q+t(1-q)] u(t)\, -\, \frac{1}{2} \rho\, t\right)\, .
\end{equation}
This quadratic equation is easily solved giving (with the correct branch)
\begin{equation}
\label{thrdeq}
  2 u(t)\, =\, -\left(\frac{q}{2q-1}\, -\, t\right)\, +\, \left[\left(\frac{q}{2q-1}-t\right)^2 + \rho\left(\rho\, +\, \frac{2t-2q}{2q-1}\right)\right]^{1/2}\, .
\end{equation}
Equations (\ref{secdeq}) and (\ref{thrdeq}) constitute the fundamental solution
of the integral equation (\ref{inteq}).

Some plots of the profiles obtained in this way are shown in Figure 1.
We have superimposed profiles for $q$ between the values of $1/2$ and $1$
with stepsize $0.025$.

It is easily seen that this does uniquely specify a measure $\mu(dx,t)$.
Therefore, there is at most one extremal thinning-invariant state describing the 
limit of the invariant-measures for the ASEP,
and we believe that it is, in fact the unique limit.

\section{A note about quantum models}

Suppose that for each $N$, one has a quantum Hamiltonian on $N$ particles,
$H_N$, derived from an $n$-body interaction, 
which we write with its localization as $(h_n)_{i_1,\dots,i_n}$,
analogous to the situation in Section 4.
Supposing that the $N$ particles are confined to a compact set,
one can apply the main theorem to the density
\begin{equation}
  \rho_N(dx_1\otimes \cdots \otimes dx_N)\, 
  =\, |\psi_0(x_1,\dots,x_N)|^2\, dx_1\otimes \cdots \otimes dx_N\, ,
\end{equation}
where $\psi_0$ is any pure ground state.
In this case, we can define $(\rho^{(N)}_k\, :\, k\leq N)$ 
as in Section 4.
The interpretation is that $\rho^{(N)}_k$ is the 
diagonal of the $k$-particle reduced density matrix,
averaged over all $N$-choose-$k$ choices of $k$ particles. 
As before, we know that any limit of ground states 
determines a thinning-invariant sequence of measures, so that
it can be represented as in Theorem 1.

However, note that the energy expectation is not a linear functional of $\rho_N$.
Instead, if the Hamiltonian is positivity-preserving,
so that we can assume the wavefunction is nonnegative,
then the energy expectation is a linear functional of 
the ground state projection
$\ket{\sqrt{\rho_N}}\bra{\sqrt{\rho_N}}$.
Therefore, we do not know how to deduce that the extreme limit
points of the ground states are 
extreme points of the simplex of thinning-invariant measures.

Of course there are some simple examples where the unique limit is an
extreme point of the thinning-invariant simplex.
For example, the following two-body interaction may be the simplest of all:
\begin{equation}
  (h_2)_{j,k} \, =\, -\left(\frac{\partial^2}{\partial x_j^2} + \frac{\partial^2}{\partial x_k^2}\right)
  + x_j^2 + x_k^2 + (x_j-x_k)\, .
\end{equation}
Because this Hamiltonian is quadratic,
generating a quasifree evolution,
the ground state factorizes for every 
choice of system size, $N$.

Although our theorem does not apply to this model because $\mathscr{X}=\R$ is not compact, 
it is trivial to see that the conclusion is still valid with 
\begin{equation}
  \mu(dx,t)\, =\, (2/\pi)^{1/2}\, \exp\left[-\sqrt{2}\left(x+t-1/2\right)^2\right]\, dx\, .
\end{equation}

\section{Subadditivity of the pressure for oriented, mean-field models}

Although we could not use our representation theorem to calculate the pressure for
classical spin systems, as defined in Section \ref{Sec:CSS}, we can prove that the
limiting pressure exists, by a very simple argument.

Suppose that $(H_N\, :\, N \geq n)$ is a sequence of Hamiltonians
for an oriented, mean-field spin chain, 
with an $n$-body interaction, as described in Section \ref{Sec:CSS}.
The normalized partition function is defined as 
\begin{equation}
  Z_N(\beta)\, =\, \sum_{\sigma \in \{-1,+1\}^N} 2^{-N} \exp(-\beta H_N(\sigma))\, ,
\end{equation}
and the thermodynamic potential is $\Omega_N(\beta) = \log Z_N(\beta)$.
We will provide an easy proof that the sequence $(\Omega_N(\beta)\, :\, N \geq n)$ is 
subadditive
in the sense that if $N_1,N_2 \geq n$,
and $N=N_1+N_2$, then
\begin{equation}
  \Omega_{N}(\beta)\, \leq\, \Omega_{N_1}(\beta)\, +\, \Omega_{N_2}(\beta)\, .
\end{equation}
This is important because it allows one to deduce the existence of the pressure:
\begin{cor}
For any model as defined in Section \ref{Sec:CSS}, 
the pressure,
\begin{equation}
  p(\beta)\, =\, \lim_{\substack{N \to \infty \\ N \geq n}} N^{-1} \log Z_N(\beta)\, ,
\end{equation}
exists, possibly as $-\infty$.
\end{cor}
The corollary is a well-known consequence of subadditivity, so we will
not prove it.
However, we will prove that in fact the pressure is finite, as follows trivially
from the Gibbs variational principle.

Subadditivity could be proved in a number of ways.
We will prove it using the Gibbs variational principle.
The advantage of this is to give an argument which is as close as possible
to the general argument for existence of the pressure for short-ranged
models on $\Z^d$.
In particular, we hope it will be clear exactly what
relationship the present problem
has to the problem of proving existence of the pressure for short-ranged
models: namely, the present problem is easier.
(There are aspects of statistical mechanics which are harder
for long-ranged models, such as determining the correct
analogue of the Dobrushin-Lanford-Ruelle equations.)

For a reference to the Gibbs variational principle,
entropy, and pressure, as they relate to spin systems,
consult \cite{Israel}, Section II.2 or \cite{Simon}, Section III.4.
These are the references we follow.

For any $N \in \N_+$, let $\Lambda_N = \{1,\dots,N\}$,
and define $\rho^0$ on $\{+1,-1\}^N$ as the uniform distribution:
$\rho^0(\{\sigma\}) = 2^{-N}$ for all $\sigma$.
Given another probability measure $\rho$ on $\{+1,-1\}^N$,
the relative entropy with respect to $\rho^0$ is 
\begin{equation}
  S_N(\rho)\, =\, \sum_{\sigma \in \{+1,-1\}^N} g\left(\frac{\partial \rho}{\partial \rho^0}(\sigma)\right)\, \rho^0(\{\sigma\}) 
\end{equation}
where we write the Radon-Nikodym derivative
\begin{equation}
  \frac{\partial \rho}{\partial \rho^0}(\sigma)\, :=\, \frac{\rho(\{\sigma\})}{\rho^0(\{\sigma\})}\, ,
\end{equation}
and $g(x) : \R_+ \to \R$ is the continuous, concave function
$g(x) = - x \log(x)$, defined by continuity at $x=0$. 
I.e., $g(0)=0$.

In classical spin systems, strong subadditivity of the entropy is an important
and well-known fact. 
We will only use subadditivity.
For any $X \subset \Lambda_N$, one can define the restriction
$\rho \restriction X$ such that for any $\sigma_X \in \{+1,-1\}^X$,
\begin{equation}
  (\rho\restriction X)(\{\sigma_X\})\, 
  =\, \sum_{\sigma \in \{+1,-1\}^N} \mathbbm{1}[\sigma\restriction X\, =\, \sigma_X]\,
  \rho(\{\sigma\})\, .
\end{equation}
For any $k \in \{0,\dots,N\}$, let $\mathcal{P}_k(\Lambda_N)$ be the collection
of cardinality-$k$ subsets of $\Lambda_N$.
Then subadditivity, says that
\begin{multline}
  \forall k \in \{1,\dots,N-1\}\, ,\ \forall X \in \mathcal{P}_k(\Lambda_N)\, ,\\
  S_N(\rho)\, \leq\, S_k(\rho \restriction X) + S_{N-k}(\rho \restriction (\Lambda_N \setminus X))\, .
\end{multline}
See Lemma II.2.1 in \cite{Israel} or Theorem III.4.2 in \cite{Simon}.

It will also be important for us, as is proved in the references, that
the relative entropy is concave.
If $\rho_1, \rho_2 \in \mathcal{M}^+_1(\{+1,-1\}^N)$ and $0<\theta<1$,
then
\begin{equation}
  S_N(\theta\cdot \rho_1 + (1-\theta)\cdot \rho_2)\, 
  \geq\, \theta\cdot S_N(\rho_1)\, +\, (1-\theta)\cdot S_N(\rho_2)\, .
\end{equation}

Equally important is the Gibbs variational principle, 
(c.f., \cite{Israel}, II.3.1), 
which says that for any $\rho \in \mathcal{M}^+_1(\{+1,-1\}^N)$, 
one has
\begin{equation}
  \Omega_N(\beta) \geq S_N(\rho) - \beta \rho(H_N)\, ,
\end{equation}
with equality iff $\rho$ is the Gibbs distribution for $\beta H_N$.
From these facts, one can conclude that for a short-ranged Hamiltonian on 
$\Z^d$, the pressure is subadditive, modulo small errors due to
surface energies.
This is the usual way that one proves the existence of pressure for 
short-ranged models, and moreover that the pressure can be approximated
in finite-volumes.
C.f., \cite{Israel}, \cite{Ruelle} or \cite{Simon}, which all devote
chapters to proving existence of the pressure for short-ranged
classical spin systems.

To prove subadditivity of the pressure in our case, 
let us suppose that $N \geq n$
and that $\rho$ is a probability measure on $\{+1,-1\}^N$.
For any $k \in \{1,\dots,N\}$, let us define a measure $(\rho_N \restriction k)$
on $\{+1,-1\}^k$ by
\begin{equation}
(\rho_N \restriction k)(\{(\sigma_1,\dots,\sigma_k)\})\, 
  :=\, \binom{N}{k}^{-1} \sum_{X \in \mathcal{P}_k(\Lambda_N)} (\rho_N \restriction X)(\{(\sigma_1,\dots,\sigma_k)\})\, .
\end{equation}
This is a convex combination of states.
Then, as long as $k \geq n$, it is easily seen that
\begin{equation}
  (\rho_N \restriction k)(H_k)\, =\, (k/N)\, \rho_N(H_N)\, .
\end{equation}
Therefore, it is obvious that
\begin{equation}
\label{energy}
  \rho_N(H_N)\, 
  =\, (\rho_N \restriction k)(H_k)\, 
  +\, (\rho_N \restriction N-k)(H_{N-k})\, ,
\end{equation}
as long as $N \geq 2n$ and $k \in \{n,\dots,N-n\}$.

On the other hand, 
using subadditivity and concavity of the entropy,
we see that
\begin{equation}
\begin{split}
  S_N(\rho)\, 
  &\leq\,  \binom{N}{k}^{-1} \sum_{X \in \mathcal{P}_k(\Lambda_N)} \left[S_k(\rho \restriction X)\, +\, S_{N-k}(\rho \restriction X^c)\right] \\
  &=\, \sum_{X \in \mathcal{P}_k(\Lambda_N)} {\binom{N}{k}}^{-1}\, S_k(\rho \restriction X)\\
  &\qquad +\, \sum_{Y \in \mathcal{P}_{N-k}(\Lambda_N)} {\binom{N}{N-k}}^{-1}\, S_{N-k}(\rho \restriction Y)\\
  &\leq S_k\left( \sum_{X \in \mathcal{P}_k(\Lambda_N)} {\binom{N}{k}}^{-1}\,  (\rho \restriction X) \right)\\ 
  &\qquad +\, S_{N-k}\left( \sum_{Y \in \mathcal{P}_{N-k}(\Lambda_N)} {\binom{N}{N-k}}^{-1}\, (\rho \restriction Y) \right)\\
  &= S_k(\rho\restriction k)\, +\, S_{N-k}(\rho\restriction N-k)\, .
\end{split}
\end{equation}
Combining this with equation (\ref{energy}) and the Gibbs variational principle, 
leads immediately to the subadditivity which was claimed.

Moreover, using the uniform measures $\rho_N^0$
in the Gibbs principle demonstrates that for every $N \geq n$, it is true that
\begin{equation}
  \Omega_N(\beta)\, \geq\, -\beta \rho_N^0(H_N)\, =\, (N/n)[-\beta \rho_n^0(\phi_n)]\, ,
\end{equation}
because $S_N(\rho_N^0)=0$.
This shows that the pressure is bounded below, hence not $-\infty$.

This result generalizes a theorem from \cite{BCG}, where existence of the pressure
for mean-field models was proved by a more complicated technique using interpolation.
(However, it should be noted that one of the purposes of \cite{BCG} was to show that
a very important interpolation technique in the theory of spin-glasses could be
applied also to non-random spin systems.)

Finally, we would like to say that there are models where one can calculate the pressure
using ideas related to those of the present paper.
For example, for the exchangeable (in distribution), random spin system
called the Sherrington-Kirkpatrick
model, a simple argument resulted in an extended variational principle in \cite{AS2}.
We would call this {\em Aizenman's extended variational principle}, which is more
generally applicable than just to the Sherrington-Kirkpatrick model.
For example, for non-random, exchangeable spin systems, Aizenman's extended
variational principle applies whenever the Hamiltonian is a convex or 
concave function of the average magnetization.

Although for the Sherrington-Kirkptatick model, 
the Euler-Lagrange equations derived by the extended variational
principle seem too hard to solve directly, for non-random models
they are readily solved.
More specifically, what is easily determined is that, while the functional
to be optimized is nonlinear, it is homogeneous of degree zero,
and the optimizers can be chosen to be extreme points.
We plan to present a careful analysis of this point,
with further applications, in another paper.

\section*{Acknowledgements}

It is a pleasure to thank Janko Gravner for the invaluable suggestion to
look for applications, and for his suggestion that I consult reference
\cite{Liggett}.
I thank Eugene Kritchevski for an interesting discussion
about potential theory (though the application did not
make it to this paper).
This research was supported by a CRM-ISM fellowship.

\end{document}